\documentclass{article}

\usepackage{amsfonts}
\usepackage{amssymb}
\usepackage{enumerate}

\newtheorem{ttt}{Theorem}[section]
\newtheorem{llll}[ttt]{Lemma}
\newtheorem{ccc}[ttt]{Claim}
\newtheorem{sss}[ttt]{Statement}
\newtheorem{ddd}[ttt]{Definition}
\newtheorem{qqq}[ttt]{Question}
\newtheorem{cccc}[ttt]{Corollary}

\newcommand{\bt}{\begin{ttt}}
\newcommand{\bl}{\begin{llll}}
\newcommand{\bc}{\begin{ccc}}
\newcommand{\bs}{\begin{sss}}
\newcommand{\bd}{\begin{ddd} \upshape}
\newcommand{\bq}{\begin{qqq}}
\newcommand{\bcor}{\begin{cccc}}

\newcommand{\bp}{\noindent\textbf{Proof }}
\newcommand{\br}{\noindent\textbf{Remark }}

\newcommand{\et}{\end{ttt}}
\newcommand{\el}{\end{llll}}
\newcommand{\ec}{\end{ccc}}
\newcommand{\es}{\end{sss}}
\newcommand{\ed}{\end{ddd}}
\newcommand{\eq}{\end{qqq}}
\newcommand{\ecor}{\end{cccc}}

\newcommand{\ep}{\hspace{\stretch{1}}$\blacksquare$\medskip}
\newcommand{\er}{\medskip}

\newcommand{\lab}[1]{\label{#1}}

\newcommand{\NN}{\mathbb{N}}

\newcommand{\QQ}{\mathbb{Q}}
\newcommand{\RR}{\mathbb{R}}

\newcommand{\oegy}{\omega_1}
\newcommand{\al}{\alpha}
\newcommand{\be}{\beta}
\newcommand{\ga}{\gamma} 
\newcommand{\de}{\delta}
 
\newcommand{\om}{\omega}

\newcommand{\iH}{\mathcal{H}} 
\newcommand{\iR}{\mathcal{R}}
\newcommand{\iF}{\mathcal{F}}
\newcommand{\iS}{\mathcal{S}}
\newcommand{\iT}{\mathcal{T}}

\title{Linearly Ordered Families of Baire 1 Functions}

\author{M. Elekes}

\begin{document}

\maketitle 

\begin{abstract}
We consider the set of Baire 1 functions endowed with the pointwise partial
ordering and investigate the structure of the linearly ordered subsets. 
\end{abstract}

\section*{Introduction}

Any set $\iF$ of real valued functions defined on an arbitrary set $X$ is
partially ordered by the pointwise ordering, that is $f\leq g$ iff
$f(x)\leq g(x)$ for all $x\in X$. In other words put $f<g$ iff $f(x)\leq
g(x)$ for all $x\in X$ and $f(x)\not= g(x)$ for at least one $x\in X$. Our
aim will be to investigate the possible order types of the linearly ordered
(or simply `ordered' from now on) subsets of this partially ordered set,
which is the same as to characterize the ordered sets that are similar to
an ordered subset of $\iF$. Here two ordered sets are said to be similar
iff there exists an order preserving bijection between them, and such a
bijection from an ordered set onto an ordered subset of $\iF$ is often
referred to as a `representation' of the ordered set. We sometimes say that
the set is represented `on $X$'. An ordered set similar to a representable
one is also representable, so we can talk about `representable order types'
as well.

Since the functions in an ordered set are somehow `above each other', one
could think that this ordered set must be similar to a subset of the real
line. As we shall see this is far from being true.

The problem of finding long sequences in $\iF$, that is representing big
ordinals has been studied for a long time.  It was Mikl\'os Laczkovich who
posed the question how one can characterize the representable ordered sets,
particularly in the case when $X=\RR$ and $\iF$ is the set of Baire 1
functions.  What makes this problem interesting is that the corresponding
questions about continuous (that is Baire 0) and Baire $\al$ functions
($\al>1$) are completely solved. In the continuous case an ordered set is
representable iff it is similar to a subset of $\RR$ (an easy exercise),
and for $\al>1$ the question has turned out to be independent of $ZFC$,
that is the usual axioms of set theory \cite{Ko}.

The known facts about the case $\al=1$ are the followings.  The first is a
classical theorem of Kuratowski asserting that there is no increasing or
decreasing sequence of length $\oegy$ of real Baire 1 functions \cite[\S
24. III.2']{Ku}, that is $\oegy$ is not representable (in the sequel
representable will always mean representable by real Baire 1 functions).
The other is P\'eter Komj\'ath's Theorem stating that no Souslin line is
representable \cite{Ko}. (A Souslin line is a non-separable ordered set
that does not contain an uncountable family of pairwise disjoint open
intervals, that is ccc but not separable. The existence of Souslin lines is
independent of $ZFC$ \cite[Theorems 48,50]{Je}.)

The main goal of this paper is to present a few constructions of
representable ordered sets which show that Kuratowski's Theorem is `not too
far' from being a characterization. In Section 2 we prove that certain
operations result representable order types, and then in Section 3 and 4 we
show that everything is representable that can be built up by certain steps,
like forming countable products or replacing points by ordered sets.

We would also like to point out that if we restrict ourselves to the case
of characteristic functions, we arrive at the problem of families of sets
linearly ordered by inclusion.  Indeed, $\chi_A<\chi_B$ iff $A\varsubsetneqq
B$. The case of real Baire 1 functions corresponds to the problem of
representing ordered sets by ambiguous subsets of the real line. (A set is
called ambiguous iff it is $F_\sigma$ and $G_\delta$ at the same time.) It
is not hard to check that almost everything proved in this paper is valid
for this case as well, moreover, a kind of characterization of ordered sets
that are representable by ambiguous sets is given in the last section.

For a topological space $X$ the set of order types representable by real
valued Baire 1 functions is denoted by $\iR(X)$.  The set of order
types representable by ambiguous subsets is denoted by $\iR_0(X)$.

\bigskip\noindent{\textbf{Acknowledgment }  I am greatly indebted to my
advisor Professor Mik\-l\'os Lacz\-ko\-vich for his much help and advice
and for
everything I have learnt from him.} 

\section{Preliminaries}

We shall frequently use the following simple lemma. 

\bl \lab{baireone}  \ 
\begin{enumerate}[(i)]
\item Let $X$ and $Y$ be metric spaces, $f:X\rightarrow\RR$ Baire 1 and
$g:Y\rightarrow X$ continuous.  Then $f\circ g:Y\rightarrow\RR$ is Baire 1. 
\item Let $X$ be a metric space and $X_n \subset X$ $(n\in\NN)$ $F_\sigma$
sets such that $X=\bigcup_{n=1}^\infty X_n$. If $f:X\rightarrow\RR$ is
relatively Baire 1 on each $X_n$ $(n\in\NN)$ then $f$ is Baire 1. 
\end{enumerate}
\el 

Let us first consider the following question, which shall be a useful tool
in the sequel. Which Polish spaces are equivalent to the real line in the
sense that the same ordered sets can be represented on them? We shall
ignore the countable metric spaces as it is easy to see that if an order
type is representable on such a space then it is similar to a subset of the
real line. Denote by $C$ the Cantor set.

\bt \lab{scpt} $R(X)=R(C)=R(\RR)$ for any $\sigma$-compact uncountable metric
space $X$.
\et 

\bp It is obviously enough to prove the first equality. Let $X$ be compact
for the time being, then a classical theorem asserts that there exists a
continuous surjection $F:C\rightarrow X$ \cite[\S 41, VI.3a]{Ku}.  If
$\{f_\alpha:\alpha\in\Gamma\}$ is an ordered set of Baire 1 functions
defined on $X$, one can easily verify that $\{ f_\alpha\circ
F:\alpha\in\Gamma\}$ is also ordered, similar to the former ordered set as
a consequence of the surjectivity of $F$ and consists of Baire 1 functions
defined on $C$ by lemma \ref{baireone}. 

In the general case $X=\cup_{n=1}^\infty X_n$ where $X_n\subset X$ is
compact and let again be $\{f_\alpha:\alpha\in\Gamma\}$ an ordered set of
Baire 1 functions on $X$. We shall show that this set is representable on
the interval $[0,1]$ and therefore on $C$ as well, since $[0,1]$ is a
compact metric space and we can apply what we have proven in the previous
case. 

Fix a set $H_n\subset (\frac{1}{n},\frac{1}{n+1})$ for each $n\in\NN$
homeomorphic to the Cantor set and also a homeomorphism $g_n:H_n\rightarrow
C$. We can choose furthermore continuous surjections $F_n:C\rightarrow X_n$
$(n\in\NN)$ since $X_n$ is a compact metric space. Now we represent the set
in the following way. For each $\alpha\in\Gamma$ let 
\[
g_\alpha= \left\{ \begin{array}{ll} f_\alpha\circ F_n\circ g_n & \textrm{on
	$H_n$ $\ (n\in\NN)$}\\ 0 & \textrm{on
	$[0,1]\setminus\cup_{n=1}^\infty H_n$.}  \end{array}\right. 
\]
Indeed, the map $g_\alpha\mapsto f_\alpha\ (\alpha\in\Gamma)$ turns out to
be a similarity as $F_n\circ g_n$ is surjective and moreover in view of
Lemma \ref{baireone} it is straightforward to verify that $g_\alpha$ is a
Baire 1 function on $[0,1]$ for each $\alpha\in\Gamma$. 

In order to check the opposite direction let $\{f_\alpha:\alpha\in\Gamma\}$
be an ordered set of Baire 1 functions on the Cantor set. According to a
classical theorem every uncountable compact metric space contains a
subspace homeomorphic to C \cite[\S 36, V.1]{Ku}, which easily generalizes
to the case of uncountable $\sigma$-compact metric spaces since if
$X=\cup_{n=1}^\infty X_n$, $X_n$ compact, then at least one $X_n$ is
uncountable. We can therefore fix a homeomorphism $h:C\rightarrow Y\subset
X$ and for $\alpha\in\Gamma$ let 
\[
g_\alpha= \left\{ \begin{array}{ll} f_\alpha\circ h^{-1} & \textrm{on
	$Y$}\\ 0 & \textrm{on $X\setminus Y$.}  \end{array}\right. 
\]
One can easily prove in the above manner that this is an ordered set of
Baire 1 functions similar to the above one.  
\ep 

The above theorem implies the surprising fact that all the complicated ordered
sets represented in the following sections are also representable by
functions of connected graphs. 

\bcor A representable ordered set is also representable by Darboux Baire 1
functions and consequently by Baire 1 functions of connected graphs.  
\ecor 

\bp It is well-known that the graph of a Baire 1 function is connected iff
it is Darboux \cite[II.1.1]{Br}. By the previous theorem we can assume that
the set is represented on the Cantor set. It is not hard to extend the
representing functions by a common continuous function to the complement of
the Cantor set which makes the representing functions Darboux and Baire 1
by Lemma \ref{baireone}.  
\ep

Next we show that there are at most two distinct possible sets $\iR(X)$ for
all uncountable Polish spaces $X$.

\bt $R(X)=R(\RR\setminus\QQ)$ for any non-$\sigma$-compact Polish space $X$.  
\et 

\bp We apply the argument of Theorem \ref{scpt}. In one direction we use that
every Polish space is the continuous
image of the irrationals \cite[\S 36, II.1]{Ku}, while in the other
direction we apply Hurewicz's Theorem \cite[Theorem 7.10]{Ke} asserting
that every non-$\sigma$-compact Polish space contains a homeomorphic copy
of the irrationals as a closed subspace.  
\ep 

This leaves the question open whether all uncountable Polish spaces are
equivalent or not.

\bq Does $R(C)=R(\RR\setminus\QQ)$ hold?  
\eq 

\br In order to give an affirmative answer it would be enough to prove that
every ordered set of Baire 1 functions
on the irrationals can be represented by Baire 1 functions on the
reals. Indeed, on one hand every uncountable Polish space contains a subset
which is homeomorphic to the Cantor set \cite[\S 36, V.1]{Ku}, and on the
other hand every Polish space is the continuous image of $\RR\setminus\QQ$
hence the above argument works. 

Moreover, it can be shown that a Baire 1 function defined on the irrationals
can be extended to the reals as a Baire 1 function, but so far we were
unable to do this in an order preserving way.  
\er 

\section{Operations on representable ordered sets}

Now we investigate whether the class of representable sets are closed under
certain operations. We shall make use of these operations when constructing
complicated representable ordered sets. 

\bd For an arbitrary ordered set $X$ we call $X \times \{0,1\}$ with the
lexicographical ordering the \emph{duplication of $X$}.  
\ed 

\bq Is it true that the duplication of a representable set is also
representable?  
\eq 

In most cases this question can be replaced by the following statement. 

\bs \lab{dupl} Let $X$ be an ordered set such that the duplication of $X$
is representable. Then so is the ordered set obtained by replacing every $x
\in X$ by a representable set $Y_x$, that is $ \{(x,y):x\in X,y \in Y_x\}$
with the lexicographical ordering.  
\es 

\bp First we replace the points of the real line by uncountable closed sets
in the following way.  Let $P:[0,1] \rightarrow [0,1]^2$ be a Peano curve,
that is a continuous surjection, and let $P_1$ be its first coordinate
function.  Then $P_1 :[0,1] \rightarrow [0,1]$ is also a continuous
surjection, moreover the preimages $P^{-1}_1(\{c\})$ are uncountable closed
sets for all $c\in[0,1]$. In virtue of Theorem \ref{scpt} we may assume
that the duplication of $X$ is represented on $[0,1]$ by the pairs of
functions $f_x<g_x$ $(x \in X)$. If we consider the functions $f_x \circ
P_1$ and $g_x \circ P_1$ we obtain a similar ordered set of Baire 1
functions, but in the latter set any two distinct elements differ on an
uncountable closed sets, for if $f_x$ and $g_x$ attained different values
at $c_x$ then $f_x \circ P_1$ and $g_x \circ P_1$ differ on
$P^{-1}_1(\{c_x\})$.  Since this is a compact metric space we may assume
that $Y_x$ is represented on it.  By composing with a increasing
homeomorphism between $\RR$ and the interval $(f_x(c_x),g_x(c_x))$ we also
can assume that the functions representing $Y_x$ only attain values between
$f_x(c_x)$ and $g_x(c_x)$. 

Now we claim that the following representation will do. For $x \in X$ and
$y \in Y_x$ let 
\[
h_{(x,y)}=\left\{ \begin{array}{ll} f_x\circ P_1 & \textrm{on
	$[0,1]\setminus P^{-1}_1(\{c_x\})$}\\ \textrm{the function
	representing y}&\textrm{on $P^{-1}_1(\{c_x\})$.}  \end{array}
	\right. 
\]
These functions are easily seen to be Baire 1 so what remains to show is
that the representation is order preserving. In the first case $x_1<x_2$ so
$f_{x_1} < g_{x_2}$ hence 
\[
h_{(x_1,y_1)} < g_{x_1}\circ P_1 < f_{x_2}\circ P_1 < h_{(x_2,y_2)}. 
\]
Finally, in the second case $x_1=x_2=x$ and $y_1<y_2$.  Obviously
$h_{(x,y_1)}$ and $h_{(x,y_2)}$ differ on $P^{-1}_1(\{c_x\})$ only, where
they are defined according to the ordering of $Y_x$ thus
$h_{(x,y_1)}<h_{(x,y_2)}$.  
\ep 

\bs \lab{product} Let $X$ be an ordered set such that the duplication of $X$
is representable. Then $X^\om$ endowed with the lexicographical ordering is
also representable.
\es

\bp
As in the previous proof we can represent the duplication of $X$ such that
for every $x\in X$ the representing functions $f_x, \ g_x : \RR \rightarrow
[0,1]$ are different constant functions on a suitable
Cantor set $C_x$. Denote $d_x$ the difference of these two values. In the
next step, for every fixed $x_1\in X$ let us represent the duplication of
$X$ on $C_{x_1}$ in the same manner as above, that is for each $x_2\in X$
let $f_{x_1,x_2}, \ g_{x_1,x_2} : \RR \rightarrow
[0,\min(\frac{1}{2},d_{x_1})]$ be zero outside $C_{x_1}$ such that they
are different constants on a suitable Cantor set $C_{x_1,x_2}\subset
C_{x_1}$. Let $d_{x_1,x_2}$ denote the difference of the two values. Then
we proceed inductively and make sure that $0\leq f_{x_1,\dots,x_{n+1}},\
g_{x_1,\dots,x_{n+1}}\leq\min(\frac{1}{2^n},d_{x_1,\dots,x_n})$. It is
not hard to see that 
\[
(x_1,x_2,\dots) \mapsto \sum^\infty_{n=1} f_{x_1,\dots,x_{n}}
\]
is the required representation, as the uniform limit of Baire 1 functions
is Baire 1 itself \cite[\S 31, VIII.2]{Ku}.
\ep

\br
Instead of using the same set $X$ at each level, we can prove in exactly
the same way that if the duplication of $X_n$ is representable for every
$n\in\NN$ then so is $\prod^\infty_{n=1} X_n$, and more generally we can
also use different sets at a level, that is we can correspond a set
$X_{x_1,\dots,x_n}$ to each $x_1,\dots,x_n$. 
\er

However,  we do not know the answer to the question concerning longer
products. As a simple transfinite induction shows, the following two questions
are equivalent.

\bq
Is it true, that if the duplication of $X$ is representable, then the
duplication of $X^\om$ is also representable? Or equivalently, is it true,
that if the duplication of $X$ is representable, then so is $X^\al$ for
every $\al<\oegy$? 
\eq

\bcor Suppose that the duplications of representable orderings are also
representable. Then $X^\alpha$ is representable for every representable $X$
and $\alpha<\omega_{1}$.
\ecor

\bp
We prove this by induction on $\al$. If $\al=\be+1$ then $X^\al$ is similar
to $X^\be\times X$. But $X^\be$ is representable by the inductional
hypothesis, so is its duplication by our assumption, therefore we can apply
Statement \ref{dupl} and we are done.

If $\al$ is a limit ordinal, then $[0,\al)$ can be written as the disjoint
union of $[\al_n,\al_{n+1})$ for a suitable sequence
$\al_n\ (n\in\NN)$. The interval $[\al_n,\al_{n+1})$ is similar to an ordinal
$\be_n<\al$, so $X^\al$ is similar to $\prod^\infty_{n=1} X^{\be_n}$, and
we are again done by the previous remark.
\ep

\br
As above, we can generalize this result as well to $\prod_{\be<\al} X^\be$
and also to the case when at each level we correspond an arbitrary
representable set to each point.
\er

Next we pose another question. 

\bq Is it true that the completion (as an ordered set) of a representable
ordered set is also representable?  
\eq 

\bd Let $X$ and $X_n\ (n\in \NN)$ be ordered sets.  We say that $X$ is
\emph{a blend of the sets $X_n$} if there exist pairwise disjoint subsets
$H_n \subset X \ (n\in\NN)$ such that $X=\cup_{n=1}^\infty H_n$ and $H_n$
is similar to $X_n$.  
\ed 

\bs Suppose that duplications and completions of representable sets are
also representable. Then so is a blend $X$ of the representable sets $X_n$.
\es 

\bp Let $H_n$ be as in the definition.  By the hypothesis the completion of
 $H_n\times \{0,1\}$ is representable for each $n\in\NN$ and we may assume
 that it is represented on the interval $(n,n+1)$.  Let $x\in X$, that is
 $x\in H_n$ for exactly one $n$, and let 
\[
f_x=\left\{\begin{array}{ll} \textrm{the function representing $(x,0)$} &
	\textrm{on $(n,n+1)$}\\ \textrm{the function representing}\\
	\textrm{$\sup\{(y,i)\in H_m\times \{0,1\}:y\leq x\}$} & \textrm{on
	$(m,m+1)\ if\ m\not=n$}\\ $0$ & \textrm{elsewhere},
	\end{array}\right. 
\]
where `sup' means supremum according to the ordering of the completion of
$H_m\times \{0,1\}$. $f_x$ is Baire 1 as the usual argument shows so we
only have to check that this latter set of functions is similar to the
original one. Let $x,y\in X$, $x<y$ and $x\in H_k$, $y\in H_l$ for some $k$
and $l$. If $k=l$ then $f_x < f_y$ is obvious while if $k\not= l$ then one
can easily check that $f_x\leq f_y$ on $(k,k+1)$, $(l,l+1)$ and on the
complement of their union, moreover $f_x\not=f_y$ on $(k,k+1)$ since $f_y$
is not less here then the function representing $(x,1)$.  
\ep 

\section{The first construction}

In the sequel we present a few constructions of representable sets which
have such a rich structure in some sense that we may hope to be able to
produce all the representable order types this way. 

\bd Let $\alpha$ be an ordinal number and $I=[0,1]$.  We denote by
$I^{\alpha}$ the set of transfinite sequences in $I$ of length $\alpha$
with the lexicographical ordering (i.e. $I^{\alpha}=\{f:\ f:\alpha\to I\}$
and $f<g$ iff $f(\gamma)=g(\gamma)$ and $f(\beta)<g(\beta)$ for some
$\beta$ and every $\gamma<\beta$).  
\ed 

When $\alpha\geq\oegy$, then due to Kuratowski's Theorem \cite[\S 24,
III.2']{Ku}, $I^\alpha$ is not representable as it contains a subset of type
$\oegy$. However the following holds. 

\bt \lab{ialpha} $I^\alpha$ is representable for all $\alpha<\oegy$.   
\et 

\bp For $\alpha < \omega$ the assertion follows from Statement \ref{dupl} by
induction. Denote by $H=\prod_{n=0}^\infty [0,1]$ the Hilbert cube, that is
the topological product of countably many copies of the closed unit
interval. It is well-known that $H$ is a compact metric space so it is
sufficient to represent $I^\alpha$ on $H$.  We show that this is possible
even by characteristic functions, in other words there exists a system of
ambiguous subsets of $H$ which is of order type $I^\alpha$ when ordered by
inclusion.  First we define an ordering of type $I^\alpha$ on $H$.  As
$\alpha < \omega _1$ there exists a bijection $\varphi:\NN \rightarrow
\alpha$ so we can assign to each element $a=(a_1,a_2,\ldots) \in H$ a
transfinite sequence $x=(a_{\varphi(n)}:n \in \NN)$.  Since this is a
bijection between $H$ and $I^\alpha$ it induces an ordering of type
$I^\alpha$ on $H$ which we shall denote by $<_H$. We claim that the sets of
the form $H_x= \{y \in H : y<_H x\}$ constitute a system of sets possessing
all the properties we need. First of all $H_x \subsetneqq H_y$ iff $x<_H y$
thus $\{H_x:x\in H\}$ is of order type $I^\alpha$.  We still have to check
that $H_x\subset H$ is ambiguous for all $x\in H$. First we show that it is
$F_\sigma$. Indeed, 
\[
H_x=\bigcup_{\beta<\alpha} \left( \bigcap_{\gamma<\beta} \left\{
(y_1,y_2,\ldots)\in H: y_{\varphi^{-1}(\gamma)}= x_{\varphi^{-1}(\gamma)}
\right\} \cap \left\{ y_{\varphi^{-1}(\beta)}< x_{\varphi^{-1}(\beta)}
\right\} \right) 
\]
so it is sufficient to check that the members of the union are $F_\sigma$
sets, but this is obvious as they are intersections of certain closed sets
and an open set. 

Similarly $\{y\in H: x<_H y\}$ is also $F_\sigma$, and as $\{x\}$ is
$F_\sigma$, $H_x$ is the complement of an $F_\sigma$ set hence $G_\delta$.
\ep 

In view of Kuratowski's Theorem it is natural to ask whether every
representable set can be embedded into $I^\alpha$ for a suitable
$\alpha<\oegy$. We show in two steps that this is not true. 

\bl $I^{\alpha+1}$ cannot be embedded into $I^\alpha$ for any
$\alpha<\oegy$.  
\el 

\bp Suppose indirectly that $f:I^{\alpha+1} \rightarrow I^\alpha$ is an
order-preserving injection and let $f=(f_0,f_1,\ldots,f_\beta,\ldots)$
where $f_\beta:I^{\alpha+1}\to I\ (\beta<\alpha)$ are the coordinate
functions.  As $f_0:I^{\alpha+1} \rightarrow I$ is monotone, and for
distinct values of $c\in I$ the convex hulls of the sets $f_0 (\{x_0,\ldots
,x_\beta,\ldots ,x_\alpha :x_0=c\})$ are non-overlapping intervals in $I$,
all but countably many of them are singletons. Therefore we can fix $a_0$
such that $f_0 ((a_0,x_1,\ldots ,x_\beta,\ldots ,x_\alpha))$ is constant.
Once we have already chosen $a_\gamma$ for each $\gamma<\beta$ such that
$f_\gamma((a_0,\ldots,a_\gamma,x_{\gamma+1},\ldots,x_\alpha))$ is constant
then as before for distinct values of $x_\beta$ we obtain essentially
pairwise disjoint image sets and thus we can fix $a_\beta\in I$ such that
$f_\beta((a_0,\ldots,a_\beta,x_{\beta+1},\ldots,x_\alpha))$ is
constant. But then eventually we get 
\[
f\left( (a_0,\ldots,a_\beta,\ldots,0) \right) = f\left(
(a_0,\ldots,a_\beta,\ldots,1) \right), 
\]
contradicting the injectivity of $f$.  
\ep 

\bs \lab{nemialfa} There exists a representable set that is not embeddable
into $I^\alpha$ for any $\alpha<\oegy$.  
\es 

\bp The duplication of the real line is representable as it is similar to a
subset of $I^2$, hence if we replace $\aleph_1$ arbitrary points of $\RR$
by the sets $I^\alpha\ (\alpha<\oegy)$ we obtain a representable set. In
virtue of the previous lemma and Statement \ref{dupl} this set possesses the
required property.  
\ep 

This negative result shows how to go on to find new representable sets
by iteration. 

\bd Let $\mathcal{H}$ be an arbitrary set of ordered sets.  We define an
increasing transfinite sequence $S_\alpha\ (\alpha\in On)$ of sets as
follows. 

Let $S_0= \mathcal{H}\cup\{\emptyset\}$ and $S_\alpha$ be the set of
ordered sets that can be obtained by replacing the points of a set $X\in
\bigcup_{\beta<\alpha}S_\beta$ by sets $Y_x\in
\bigcup_{\beta<\alpha}S_\beta\ (x\in X)$. 

Finally, let $\mathcal{S}(\mathcal{H})$ denote the set of order types of
$\bigcup_{\alpha\in On} S_\alpha$.  
\ed 

\bl $\mathcal{S}(\mathcal{H})$ is a set indeed as there exists an ordinal
$\alpha$ such that $S_\beta=S_\alpha$ for every $\beta\ge\alpha$.  
\el 

\bp Let $\kappa$ be a infinite cardinal such that $|H|\le \kappa$ for every
$H\in\mathcal{H}$.  A simple transfinite induction shows that $|X| \le
\kappa$ for all $X\in S_\alpha$ and $\alpha\in On$.  We choose a cardinal
$\mu$ of cofinality greater than $\kappa$ (e.g.  $2^\kappa$), and claim that
$\alpha=\mu$ will do. 

 First we show that $S_\alpha = \bigcup_{\beta<\alpha} S_\beta$.  Choose
$X\in S_\alpha$, that is $Y, Z_y\in \bigcup_{\beta<\alpha}S_\beta$ and fix
$\beta, \beta_y<\alpha\ (y\in Y)$ such that $Y\in S_\beta$ and $Z_y\in
S_{\beta_y}\ (y\in Y)$.  The set $\{\beta\}\cup\{\beta_y:y\in Y\}$ is at
most of power $\kappa$ which is less then the cofinality of $\alpha$ thus
we can find a $\beta^*<\alpha$ such that $\beta, \beta_y<\beta^*\ (y\in
Y)$. But then $X\in S_{\beta^*}\subset \bigcup_{\beta<\alpha} S_\beta$. 

Secondly, we check by transfinite induction that $S_\be =S_\al$ for all
$\be \ge \al$. Suppose $S_\ga =S_\al$ for $\al\le\ga<\be$ and let $X\in
S_\be$, that is $Y, Z_y\in \bigcup_{\ga<\be} S_\ga$. However, 
\[
\bigcup_{\ga<\be} S_\ga = \bigcup_{\ga<\be} S_\al = S_\al =
\bigcup_{\de<\al} S_\de 
\] 
which implies $X\in S_\al$ by repeating the above argument.  
\ep 

\bt \lab{salpha} If $\iH$ is a set of ordered sets such that the
duplications of the elements of $\iH$ are representable, then the elements
of $\iS(\iH)$ are also representable.  
\et 

\bp We prove by transfinite induction on $\al$ the seemingly stronger
statement that even the duplications of elements of $\iS(\iH)$ are
representable.  For $\al=0$ this is just a reformulation of our
assumption. Suppose now that the statement holds for all $\be<\al$ and let
$X\in S_\al$, that is $Y, Z_y\in\bigcup_{\beta<\alpha}S_\beta$.  As
$Z_y\in\bigcup_{\beta<\alpha}S_\beta$ $Z_y\times \{0,1\}$ is representable
by the inductional hypothesis. Moreover if we replace the points of $Y$ by
the sets $Z_y\times \{0,1\}$ what we obtain is exactly the duplication of
$X$, which therefore turns out to be representable as by the inductional
hypothesis $Y\times \{0,1\}$ is representable and so we can apply Statement
\ref{dupl}.  
\ep 

\bd
If $\iH$ is a set of ordered sets, then let
\[\iH^\om=\{Y: Y\subset X^\om, X\in\iH\},
\]
and let $\iH^*$ be the closure of $\iH$ under the operations $X\mapsto
X^\al\ (\al<\oegy)$. (This closure can be formed by a similar transfinite
construction as $\iS(\iH)$.)
\ed

\bcor
 If $\iH$ is a set of ordered sets such that the
duplications of the elements of $\iH$ are representable, then the elements
of $\iS(\iH)^\om$ are also representable. This holds even for $\iS(\iH)^*$,
 assuming that the duplications of representable sets are representable.
\ecor

\br (a) We could define similar notions with products instead of powers, or
even with the more complex constructions mentioned in the remark following
Statement \ref{product}, but in fact we would not get more, as in the case
we are interested in, there are always at most continuum many sets
involved, thus we can put them together (e.g. replace the points of $\RR$
by them) to form a huge set $X$ that contains each of them, and so
the power of this set $X$ contains subsets similar to all these above
constructions.

(b) If we begin our procedure of building large representable orderings, we can
start with some set of simple ordered sets,
for example the ones representable by constants or even continuous
functions.  In both cases we have $\iH=\{\RR\}$. It is not hard to prove
that we will not get too far this way as 
$I^\om$ will not be in $\iS(\iH)$. (The proof goes by transfinite
induction. Note that any non-trivial subinterval of $I^\om$ contains a copy
of $I^\om$ and that building up a set $X$ by replacing each element $y$ of
a set $Y$ by $X_y$ is the same as partitioning $X$ into subintervals that
are ordered similarly to $Y$ such that each subinterval is similar to the
corresponding $X_y$.)  Therefore we prefer starting with the set of
`unboundedly wide trees', $\{I^\al:\al<\oegy\}$. 

(c) According to the previous theorems $\iS(\{I^\al:\al<\oegy\})$ contains
order types of representable duplication only, as the duplication of
$I^\al$ is a subset of $I^{\al+1}$. However, $\iS(\{I^\al:\al<\oegy\})\neq
\iR(\RR)$ as every element of the former set contains a
non-trivial subinterval that is similar to a subset of $I^\al$ for some
$\al$, while if $X$ is as in the proof of Statement \ref{nemialfa}, then
$X^\om$ does not. Therefore $\iS(\{I^\al:\al<\oegy\})^\om$ is a strictly
larger class of representable orderings. This holds for
$\iS(\{I^\al:\al<\oegy\})^*$ as well, under the assumption about duplications. 

It seems quite plausible that if we are allowed to replace
points by arbitrarily large sets of the form $I^\al$ (of course
$\al<\oegy$), and allowed to form countable products, then we can build up
every set not containing a sequence of
length $\oegy$.  Moreover it can be shown that $\iS(\{I^\al:\al<\oegy\})^*$
is closed under duplication, completion and blends. (The definition of
these notions for order types instead of ordered sets is obvious.) Together 
with Kuratowski's
Theorem this motivates the following question.  
\er 

\bq Does either $\iS(\{I^\al:\al<\oegy\})^\om=\iR(\RR)$ or
$\iS(\{I^\al:\al<\oegy\})^*=\iR(\RR)$ hold?  
\eq 

\section{The second construction}

Now we turn to an other approach of the problem which results in a notion
very similar to $\iS(\iH)$. 

\bs Let $\{f_\alpha:\alpha\in\Gamma\}$ be an ordered set of functions
defined on a second countable topological space and possessing the Baire
property. If any two functions differ on a set of second category then the
ordered set is similar to a subset of the real line.  
\es 

\bp Recall that an ordered set is similar to a subset of $\RR$ iff it is
separable and does not contain more than countably many pairs of
consecutive elements. 

First we prove separability.  Let $X$ be the second countable space and
suppose for the time being that $X$ is a Baire space, that is every
non-empty open subset is of second category. Denote by $B$ a countable base of
the space not containing the empty set.  We construct a countable dense
subset $M$ of $\{f_\alpha:\alpha\in\Gamma\}$ in the following way. If for
$U,V\in B$ and $p,q\in\QQ$ there exists $h\in\{f_\alpha:\alpha\in\Gamma\}$
such that $p<h$ on a residual subset of $U$ and $h<q$ on a residual subset
of $V$ then we choose such an $h$. $M$ is obviously countable and to verify
that it is dense let $(f,g)$ be an open interval of the ordered set. If
this interval is empty then we are done so we may assume that there exists
an element $h_0$ of the ordered set in the interval. Obviously 
\[
X(f<h_0)=\bigcup_{p\in\QQ} X(f<p<h_0) 
\]
and 
\[
X(h_0<g)=\bigcup_{q\in\QQ} X(h_0<q<g), 
\]
where the sets on the left hand side are by assumption of second category
hence for some $p$ and $q\ X(f<p<h_0)$ and $X(h_0<q<g)$ are of second
category as well. It is easy to see that a set of second category which
also possesses the Baire property is residual in some non-empty open
subset, moreover this open set can be chosen to be an element of $B$. As
$f,g$ and $h_0$ have the Baire property $X(f<p<h_0)$ and $X(h_0<q<g)$ have
it as well so we can find $U,V\in B$ in which these sets are residual
respectively.  But this means that for $U,V\in B$ and $p,q\in\QQ$ there
exists an element of the ordered set, namely $h_0$, satisfying all the
conditions of the definition of $M$ so there must be such an element $h\in
M$ as well. We show that $h\in(f,g)$. X is a Baire space hence $U$ is not
of first category therefore there exists $x\in U$ for which $f(x)<p<h(x)$
and similarly $y\in V$ for which $h(y)<q<g(y)$. But this implies $f<h<g$
proving the separability. 

Let now $f_i <g_i$ $(i\in I)$ be distinct consecutive elements in the
ordered set. Like above, for every $i\in I$ 
\[
X(f_i< g_i)=\bigcup_{p\in\QQ} X(f_i<p<g_i) 
\]
hence for a suitable $p_i$ $X(f_i<p_i<g_i)$ is of second category and we
can thus fix $U_i \in B$ in which this set is residual.  We show that the
map $i\mapsto (p_i,U_i)$ is injective which implies that $I$ is countable.
Indeed, if $i\not= i'$ and $(p_i,U_i)=(p_{i'},U_{i'})=(p,U)$ than, as $U$
is of second category, we obtain that for some $x\in U$ $f_i(x)<p<g_i(x)$
and $f_{i'}(x)<p<g_{i'}(x)$ contradicting the consecutiveness of the pairs. 

Finally, if $X$ is not a Baire space than as a consequence of Banach's
Union Theorem \cite[\S 10, III]{Ku} we can write it as $X=G\cup A$ where
$G$ is an open subset which is a Baire space as a subspace and $A$ is of
first category. If we consider the restrictions of the functions to $G$ we
obtain a similar ordered set as any two functions differ on a set of second
category in $X$ hence they can not coincide on $G$. In fact, by the same
argument they differ in $G$ on a set of second category and thus we can
apply what we have proven in the previous case.
\ep 

This statement enables us to simplify the structure of a represented set $X$
in the following way. Zorn's lemma implies that we can find a maximal
subset of $X$ in which every two elements differ on a set of second
category. As this subset must be separable we can choose a countable dense
subset $M$ of it. The maximal intervals of $X\setminus M$ are of a simpler
structure than $X$ since any two elements of such an interval coincide on a
residual set, moreover it follows from Kuratowski's Theorem that all
elements of the interval coincide on a common residual set.  We can thus go
on and repeat this procedure inside this residual set. This motivates the
following. 

\bd Let $\iH$ be an arbitrary set of ordered sets.  We call elements of
$\iH$ and the empty set sets of rank 0. For an ordinal $\al$ we say that an
ordered set $X$ is of rank at most $\al$ if there exists a countable subset
$M\subset X$ such that all maximal intervals $I$ of $X\setminus M$ are of
rank at most $\be$ for some $\be<\al$ where $\be$ may depend on $I$. The
class of ordered sets of rank at most $\al$ is denoted by $T_\al$. 

Finally, let $\iT(\iH)$ be the set of order types of $\bigcup_{\al\in
On}T_\al$.
\ed 

\bl If $X$ is a set of rank at most $\al$ then it is similar to a set
obtained by
replacing the points of $\RR$ by elements of $\bigcup_{\be<\al}T_\be$.
\el 

\bp Let $M\subset X$ be the countable subset as in the definition. Recall
that every countable ordered set can be embedded into $\QQ$ and fix a
$\varphi:M\to\QQ$ order preserving injective map. 

A maximal interval $I$ of $X\setminus M$ splits $M$ into two parts $M_1$
and $M_2$ in a natural way. Define 
\[
F(I)=\sup\{\varphi(x) : x\in M_1\}, 
\]
where we may assume the supremum to be finite as we may attach a first and
a last element to $X$ which may also be elements of $M$. Now if $I_1,I_2$
and $I_3$ are distinct maximal intervals following each other in this order
then we can find an element $x\in M$ between $I_1$ and $I_2$ and $y\in M$
between $I_2$ and $I_3$ therefore $F(I_1)<F(I_3)$ as
$\varphi(x)<\varphi(y)$.  Similarly, $F(I_1)=F(I_2)$ implies that there is
exactly one $x\in M$ between $I_1$ and $I_2$.  Consequently we can map $X$
to the real line via $\varphi$ and $F$ in an order preserving way such that
the preimage of a real number is one of the followings: the empty set, a
single point, a maximal interval, a maximal interval plus an extra point to
the left or right or two intervals and a point in between. But these sets
are obviously elements of $\bigcup_{\be<\al} T_\be$ hence the lemma
follows.
\ep 

\bcor If $\RR\in\iH$ then $\iT(\iH)\subset\iS(\iH)$ thus $\iT(\iH)$ is
a set indeed.
\ecor 

\bcor If the duplication of every element of rank 0 is representable then
so is every element of $\iT(\iH)$.
\ecor 

\br $\iT(\iH)= \iS(\iH)$ fails in general as the examples $\iH=\{\RR\}$ or
$\iH=\{X:X\subset I^\om\}$ show, since in both cases
$\iT(\iH)$ is a subset of the order types of $\{X:X\subset I^\om\}$.
\er 

However, the following question is open. 

\bq Does $\iS(\{I^\al:\al<\oegy\})=\iT (\{I^\al:\al<\oegy\})$ or
$\iS(\{I^\al:\al<\oegy\})^\om=\iT (\{I^\al:\al<\oegy\})^\om$ or
$\iS(\{I^\al:\al<\oegy\})^*=\iT (\{I^\al:\al<\oegy\})^*$ hold?
\eq 

\section{Final remarks}

First we give a characterization of $\iR_0(\RR)$, which in fact does not show too much about the structure of
these orderings. This is motivated by the way our constructions worked. 

\bt An ordered set $X$ is representable by ambiguous sets iff there exists
an ordering on a compact metric space such that certain initial segments
are ambiguous and ordered similarly to $X$ by inclusion.
\et 

\bp If we have such an ordering then of course the initial segments will
do.  Conversely, let $\{H_x:x\in X\}$ be a representation by ambiguous
sets. Let 
\[
a\prec b \textrm{ iff } \exists x\in X \textrm{ such that } a\in H_x
\textrm{ and } b\notin H_x. 
\]
One can easily see that this is a partial ordering on the compact metric
space. By Zorn's lemma every partial ordering can be extended to an
ordering, thus denote $\prec^*$ such an extension.  We only have to show
that $H_x$ is an initial segment indeed of $\prec^*$ for  each $x\in X$. So
let $a\in H_x$, $b\prec^* a$ and show that $b\in H_x$. If this was not true
then $b\notin H_x$, $a\in H_x$ and $b\prec^* a$ would hold, which
contradicts the definition of $\prec^*$.
\ep 

\bq Does $\iR(\RR)=\iR_0(\RR)$ hold?
\eq

To summarize our results we may say that the class of representable ordered
sets seems to be quite close to the ones not containing sequences of
length $\oegy$. Our last theorem asserts that one actually can not prove in
$ZFC$ that these two classes coincide.

\bt
The statement that a set is representable iff it does not contain a
sequence of length $\oegy$ is not provable in $ZFC$.
\et

\bp
A Souslin line does not contain such a long increasing sequence otherwise
$\{(x_\al,x_{\al+2}):\al<\oegy \ \textrm{is a limit ordinal}\}$ would be an
uncountable system of pairwise disjoint non-empty open intervals. The case
of decreasing sequences is similar. Therefore in view of Komj\'ath's
Theorem and the independence of the existence of Souslin
lines the theorem follows. 
\ep

Finally we pose a fundamental question.

\bq
Is it consistent with $ZFC$ that an ordered set is representable iff it
does not contain a sequence of length $\oegy$?
\eq

\bigskip

DEPARTMENT OF ANALYSIS
 
LOR\'AND E\"OTV\"OS UNIVERSITY

KECSKEM\'ETI U. 10-12.

H-1053 BUDAPEST, HUNGARY

E-MAIL: EMARCI@CS.ELTE.HU

\end{document}